\documentclass[english,11pt,reqno]{smfart}

\usepackage{color}
\usepackage{amsthm}
\usepackage{amsmath}
\usepackage{amsfonts}
\usepackage{amstext}
\usepackage[latin1]{inputenc}
\usepackage{amscd}
\usepackage{latexsym}
\usepackage{bm}
\usepackage{amssymb}
\usepackage[all]{xy}
\usepackage{euscript}
\usepackage{a4wide}
\usepackage{amsmath,amssymb,graphicx}
\usepackage{amssymb}
\usepackage{amsmath}
\usepackage{mathrsfs,mathtools}

\newtheorem{theorem}{Theorem}

\newtheorem{lemma}{Lemma}

\newtheorem{definition}{Definition}
\newtheorem{remark}{Remark}

\def\di{\displaystyle}

\def\eps{\varepsilon}

\newcommand{\N}{\mathbb{N}}
\newcommand{\R}{\mathbb{R}}
\newcommand{\RN}{(\mathbb{R}^d)^{N+1}}

\newcommand{\LL}{\mathcal{L}}
\newcommand{\CC}{\mathscr{C}}
\newcommand{\Q}{\boldsymbol{Q}}
\newcommand{\W}{\boldsymbol{W}}

\newcommand{\T}{\boldsymbol{T}}

\newcommand{\DM}{D^\alpha_-}
\newcommand{\DP}{D^\alpha_+}

\newcommand{\DDM}{\Delta^\alpha_-}
\newcommand{\DDP}{\Delta^\alpha_+}

\newcommand{\fonction}[5]{\begin{array}[t]{lrcl}#1 :&#2 &\longrightarrow &#3\\&#4& \longmapsto &#5 \end{array}}

\newcommand{\fonctionsansdef}[3]{\begin{array}[t]{lrcl}#1 :&#2 &\longrightarrow &#3 \end{array}}

\begin{document}
\title{Variational integrator for fractional Euler-Lagrange equations.}
\author{Lo\"ic Bourdin}
\address{Laboratoire de Math\'ematiques et de leurs Applications - Pau (LMAP). UMR CNRS 5142. Universit\'e de Pau et des Pays de l'Adour.}
\email{bourdin.l@etud.univ-pau.fr}

\author{Jacky Cresson}
\address{Laboratoire de Math\'ematiques et de leurs Applications - Pau (LMAP). UMR CNRS 5142. Universit\'e de Pau et des Pays de l'Adour.}
\email{jacky.cresson@univ-pau.fr}

\author{Isabelle Greff}
\address{Laboratoire de Math\'ematiques et de leurs Applications - Pau (LMAP). UMR CNRS 5142. Universit\'e de Pau et des Pays de l'Adour.}
\email{isabelle.greff@univ-pau.fr}

\author{Pierre Inizan}
\address{Institut de M\'ecanique C\'eleste et de Calcul des \'eph\'em\'erides, Observatoire de Paris, 77 avenue Denfert-Rochereau, 75014 Paris, France.}
\email{inizan@imcce.fr}

\date{}
\maketitle.

\begin{abstract}
We extend the notion of variational integrator for classical Euler-Lagrange equations to the fractional ones. As in the classical case, we prove that the variational integrator allows to preserve Noether-type results at the discrete level.
\end{abstract}



\section{Introduction}\label{section1}
Fractional calculus is the emerging mathematical field dealing with the generalization of the derivative to any real order. During the last two decades, it has been successfully applied to problems in economics \cite{comt}, computational biology \cite{magi} and several fields in Physics \cite{alme,bagl,hilf3,neel}. We refer to \cite{kilb,podl,samk} for a general theory and to \cite{mach} for more details concerning the recent history of fractional calculus. Particularly, a subtopic of the fractional calculus has recently gained importance: it concerns the variational principles on functionals involving fractional derivatives. This leads to the statement of fractional Euler-Lagrange equations, see \cite{agra,bale2,riew}. \\

Fractional Euler-Lagrange equations are difficult to solve explicitly and consequently, it is of interest to develop efficient numerical schemes for such dynamical systems. There exists a suitable method for classical Euler-Lagrange equations called \textit{variational integrator} and well-developed in \cite{lubi,mars}. The basic idea is to preserve at the discrete level the intrinsic variational structure of the differential equation. In this paper, we extend the notion of variational integrator to the fractional case. \\

In the classical case, a variational integrator presents many interesting numerical properties as the conservation at the discrete level of constants of motion given by Noether's theorem, see \cite{lubi,mars} for more details. In this paper, we prove that the variational integrator extended to the fractional case allows to preserve a fractional Noether-type result proved in \cite{cres6,torr3}. \\

The paper is organized as follows. Section~\ref{section2} is devoted to a reminder concerning fractional calculus and fractional Euler-Lagrange equations. In Section~\ref{section3}, we extend the notion of variational integrator to the fractional case and provide numerical simulations for the fractional Dirichlet example. In Section~\ref{section4}, we remind the fractional Noether-type result proved in \cite{cres6,torr3} and we prove that it is preserved at the discrete level by the variational integrator constructed in the previous section. We conclude this section with the fractional harmonic oscillator example. 

\section{Reminder about fractional Euler-Lagrange equations}\label{section2}
In this paper, we consider fractional differential systems in $\R^d$ where $d \in \N^*$ is the dimension. The trajectories of these systems are curves $q \in \CC^0 ([a,b],\R^d)$ with $a < b$ two reals. Let $0 < \alpha \leq 1$ denote a fractional order and let us define $\alpha_0 := 1$ and for any $r \in \N^*$:
\begin{equation}
\alpha_r := \dfrac{(-\alpha)(1-\alpha) \ldots (r-1-\alpha)}{r!}.
\end{equation}

\subsection{Fractional derivatives of Gr\"unwald-Letnikov}\label{section21}
Fractional calculus deals with the generalization of the usual notion of derivative to any real order. Since 1695, numerous notions of fractional derivatives emerge over the year, see \cite{kilb,podl,samk}. We will use the fractional operators of Gr\"unwald-Letnikov (1867) whose definitions are reminded below. \\
 
Let $q$ be a function defined on $[a,b]$ with values in $\R^d$. The left (resp. right) fractional derivative in the sense of Gr\"unwald-Letnikov with inferior limit $a$ (resp. superior limit $b$) of order $ \alpha $ of $q$ is given by:
\begin{equation}
\forall t \in ]a,b], \; \DM q (t) := \lim\limits_{\substack{h \to 0^+ \\ nh = t-a}} \dfrac{1}{h^{\alpha}} \di \sum_{r=0}^{n} \alpha_r q(t-rh),
\end{equation}
respectively:
\begin{equation}
\forall t \in [a,b[, \; \DP q (t) := \lim\limits_{\substack{h \to 0^+ \\ nh = b-t}} \dfrac{1}{h^{\alpha}} \di \sum_{r=0}^{n} \alpha_r q(t+rh),
\end{equation}
provided the right side terms are defined. \\

In the classical case $\alpha =1$, we note that $\alpha_r = 0$ for any $r \geq 2$. Consequently, we obtain that $D^1_-$ (resp. $-D^1_+$) coincides with the usual notion of left (resp. right) derivative of a function. \\

We remind the following result used further in the paper:
\begin{lemma}[Fractional integration by parts]\label{lem1}
For any $q$, $w \in \CC^1 ([a,b],\R^d)$ satisfying $w(a) = w(b) = 0$, the following equality holds:
\begin{equation}\label{fibp}\tag{FIBP}
\di \int_a^b q \cdot \DM w \; dt = \di \int_a^b \DP q \cdot w \; dt.
\end{equation}
\end{lemma}
Let us recall that the fractional derivatives of Gr\"unwald-Letnikov coincide with the Riemann-Liouville's ones (1832) for smooth enough functions, see \cite{podl}. Consequently, we refer to \cite[p.46]{samk} for a detailed proof of Lemma~\ref{lem1}. 

\subsection{Fractional Euler-Lagrange equations}\label{section22}
Since 1997 (see \cite{riew}), a subtopic of fractional calculus gains importance: it concerns the calculus of variations on functionals involving fractional derivatives. It leads to the statement of fractional Euler-Lagrange equations that was developed in \cite{agra}. \\

Let $L$ be a \textit{Lagrangian}, \textit{i.e.} a $\CC^2$ map of the form:
\begin{equation}
\fonction{L}{\R^d \times \R^d \times [a,b]}{\R}{(x,v,t)}{L(x,v,t)}
\end{equation}
and let us consider a \textit{fractional Lagrangian functional} given by:
\begin{equation}
\fonction{\LL^\alpha}{\CC^2 ([a,b],\R^d )}{\R}{q}{\di \int_a^b L(q,\DM q,t) \; dt .}
\end{equation}
Let $\CC^2_0 ([a,b],\R^d) := \{ w \in \CC^2 ([a,b],\R^d ), \; w(a) = w(b) = 0 \}$ denote the set of \textit{variations} of $\LL^\alpha$. The \textit{critical points} of $\LL^\alpha$ are the elements $q \in \CC^2 ([a,b],\R^d )$ satisfying $D\LL^\alpha (q)(w) = 0$ for any variations $w$ where:
\begin{equation}
D\LL^\alpha (q)(w) := \lim\limits_{\eps \to 0} \dfrac{\LL^\alpha (q+\eps w) - \LL^\alpha (q)}{\eps}.
\end{equation}
Finally, combining Lemma~\ref{lem1} and an usual calculus of variations, critical points of $\LL^\alpha$ are characterized by solutions of the following \textit{fractional Euler-Lagrange equation}:
\begin{equation}\tag{EL${}^\alpha$}\label{elf}
\dfrac{\partial L}{\partial x} (q,\DM q,t) + \DP \left( \dfrac{\partial L}{\partial v} (q,\DM q,t) \right) = 0.
\end{equation}
We refer to \cite{agra} for a detailed proof.

\section{Discrete fractional Euler-Lagrange equations}\label{section3}
In this section, we are interested in the construction of a variational integrator for fractional Euler-Lagrange equations \eqref{elf}. Then, let us introduce some notations available in the whole paper. Let $N \geq 2$, let $h := (b-a) / N$ be the step size of the discretization and let $\T := (t_k)_{k=0,\ldots,N} := (a+kh)_{k=0,\ldots,N}$ be the usual partition of the interval $[a,b]$.

\subsection{Discrete fractional derivatives of Gr\"unwald-Letnikov}\label{section31}
In order to give a discrete version of a fractional differential equation, we need discrete fractional operators approximating the continuous ones. As it is done in \cite{dubo}, we define the following left (resp. right) discrete fractional derivative of Gr\"unwald-Letnikov:
\begin{equation}
\fonction{\DDM}{(\R^d)^{N+1}}{(\R^d)^{N}}{\Q = (Q_k)_{k=0,\ldots,N}}{\left( \dfrac{1}{h^\alpha} \di \sum_{r=0}^k \alpha_r Q_{k-r} \right)_{k=1,\ldots,N},}
\end{equation}
respectively:
\begin{equation}
\fonction{\DDP}{(\R^d)^{N+1}}{(\R^d)^{N}}{\Q = (Q_k)_{k=0,\ldots,N}}{\left( \dfrac{1}{h^\alpha} \di \sum_{r=0}^{N-k} \alpha_r Q_{k+r} \right)_{k=0,\ldots,N-1}.}
\end{equation}
$\DDM$ (resp. $\DDP$) is an approximation of $\DM$ (resp. $\DP$). We refer to \cite{diet,lubi} for details about these approximations of fractional operators. \\

In the classical case $\alpha =1$, we note that $\Delta^1_-$ and $-\Delta^1_+$ coincide respectively with the backward and forward Euler's approximations of the derivative $d/dt$. Indeed, we have:
\begin{equation}
\fonction{\Delta^1_-}{(\R^d)^{N+1}}{(\R^d)^{N}}{\Q = (Q_k)_{k=0,\ldots,N}}{\left(  \dfrac{Q_k - Q_{k-1}}{h} \right)_{k=1,\ldots,N},}
\end{equation}
and:
\begin{equation}
\fonction{\Delta^1_+}{(\R^d)^{N+1}}{(\R^d)^{N}}{\Q = (Q_k)_{k=0,\ldots,N}}{\left( \dfrac{Q_k - Q_{k+1}}{h}  \right)_{k=0,\ldots,N-1}.}
\end{equation}

\subsection{Remark about direct discretization and variational integrator}\label{section32}
An usual and algebraic way to obtain a discrete version of a fractional differential equation is to replace the curves $q \in \CC^0 ([a,b],\R^d)$ by discrete elements $\Q \in (\R^d)^{N+1}$ and to replace fractional derivatives $D^\alpha_{\pm}$ by discrete fractional derivatives $\Delta^\alpha_{\pm}$. This method is widely used and one can find examples in \cite{agra4,liu,meer} for fractional differential equations and fractional partial differential equations. Other examples can be found in \cite{agra2,agra3,deft} for fractional optimal control problems. \\

Such a discretization leads to the following numerical scheme for a fractional Euler-Lagrange equation \eqref{elf}:
\begin{equation}\label{elfhpasbon}
\forall k=1,\ldots,N-1, \; \dfrac{\partial L}{\partial x} \big( Q_k,(\DDM \Q)_k,t_k \big) + \DDP \left( \dfrac{\partial L}{\partial v} ( \Q,\DDM \Q, \T ) \right)_k = 0.
\end{equation}
However, a fractional Euler-Lagrange equation \eqref{elf} admits a variational structure in the sense that it derives from a variational principle on a functional, see Section~\ref{section22}. This structure is intrinsic and induces strong constraints on the qualitative behaviour of the solutions. It seems then important to preserve this structure at the discrete level. Nevertheless, the numerical scheme \eqref{elfhpasbon} is obtained by an \textit{algebraic} procedure only based on the differential writing of \eqref{elf}. Consequently, there are no guarantees that the intrinsic variational structure of \eqref{elf} is preserved at the discrete level. \\

In Section \ref{section33}, we are going to construct a \textit{variational integrator} for fractional Euler-Lagrange equations \eqref{elf}. This procedure has a variational approach and allows to preserve the variational structure of \eqref{elf} at the discrete level. \\

Let us give more details concerning the construction of a variational integrator. Let us consider a differential system admitting a variational structure (\textit{i.e.} deriving from a variational principle on a functional and then, characterizing its critical points). A variational integrator is a numerical scheme constructed as follows:
\begin{itemize}
\item the first step consists in defining a discrete version of the functional;
\item the second step consists in forming a discrete variational principle on the discrete functional characterizing its discrete critical points.
\end{itemize}
Hence, a numerical scheme is obtained and it is called variational integrator. The variational structure of the differential system is then preserved at the discrete level in the sense that the discrete solutions correspond to the discrete critical points of the discrete version of the initial functional. \\

We recall that a variational integrator is well-studied for classical Euler-Lagrange equations ($\alpha =1$) in \cite{lubi,mars}. In this case, the conservation of the variational structure at the discrete level allows to preserve properties and results relative to this structure. As an example, we can cite Noether's theorem. In the following section, our aim is to extend this variational integrator to the fractional case ($0 < \alpha \leq 1$).

\subsection{Construction of a variational integrator}\label{section33}
In this section, we construct a variational integrator for fractional Euler-Lagrange equations following the two steps described in the previous section.

\paragraph*{\textbf{First step}} Let us define a discrete version of the fractional Lagrangian functional $\LL^\alpha$. Considering the usual Gaussian quadrature formula, we define the following \textit{discrete fractional Lagrangian functional}:
\begin{equation}
\fonction{\LL^\alpha_h}{(\R^d)^{N+1}}{\R}{\Q}{ h \di \sum_{k=1}^N L \big( Q_k, (\DDM \Q)_k ,t_k \big).}
\end{equation}
Hence, the first step of construction of a variational integrator for \eqref{elf} is completed. \\

\paragraph*{\textbf{Second step}} Let $(\R^d)^{N+1}_0 : = \{ \W \in (\R^d)^{N+1}, \; W_0 = W_N = 0 \}$ denote the set of \textit{discrete variations} of $\LL^\alpha_h$. In the sequel, we focus on \textit{discrete critical points} of $\LL^\alpha_h$, \textit{i.e.} elements  $\Q \in (\R^d)^{N+1}$ satisfying $D\LL^\alpha_h (\Q)(\W) = 0$ for any discrete variations $\W$ where:
\begin{equation}
D\LL^\alpha_h (\Q)(\W) := \lim\limits_{\eps \to 0} \dfrac{\LL^\alpha_h (\Q+\eps \W) - \LL^\alpha_h (\Q)}{\eps}.
\end{equation}
As in the continuous case, in order to obtain a characterization of the discrete critical points of $\LL^\alpha_h$, we first need a preliminary result. We prove the following discrete version of Lemma~\ref{lem1}:

\begin{lemma}[Discrete fractional integration by parts]\label{lem2}
For any $\Q$, $\W \in (\R^d)^{N+1}$ satisfying $W_0 = W_N = 0$, the following equality holds:
\begin{equation}\label{dfibp}\tag{DFIBP}
h \di \sum_{k=1}^N Q_k \cdot (\DDM \W)_k = h \di \sum_{k=0}^{N-1} (\DDP \Q)_k \cdot W_k.
\end{equation}
\end{lemma}

\begin{proof}
Since $W_0 = W_N = 0$, the following equalities hold:
\begin{equation}
\begin{array}{ccccc}
h^\alpha \di \sum_{k=1}^N Q_k \cdot (\DDM \W)_k & = & \di \sum_{k=1}^N \sum_{r=0}^k \alpha_r Q_k \cdot W_{k-r} & = & \di \sum_{k=0}^N \sum_{r=0}^k \alpha_r Q_k \cdot W_{k-r} \\
& = & \di \sum_{r=0}^N \sum_{k=r}^N \alpha_r Q_k \cdot W_{k-r} & = & \di \sum_{r=0}^N \sum_{k=0}^{N-r} \alpha_r Q_{k+r} \cdot W_k \\
& = & \di \sum_{k=0}^N \sum_{r=0}^{N-k} \alpha_r Q_{k+r} \cdot W_k & = & \di \sum_{k=0}^{N-1} \di \sum_{r=0}^{N-k} \alpha_r Q_{k+r} \cdot W_k \\
& & & = & h^\alpha \di \sum_{k=0}^{N-1} (\DDP \Q)_k \cdot W_k.
\end{array}
\end{equation}
Multiplying by $h^{1-\alpha}$, the proof is complete.
\end{proof}

Finally, combining Lemma~\ref{lem2} and a discrete calculus of variations, we obtain the following result:

\begin{theorem}\label{thm1}
Let $\Q \in (\R^d)^{N+1}$. Then, $\Q$ is a discrete critical point of $\LL^\alpha_h$ if and only if $\Q$ is solution of the following \textit{discrete fractional Euler-Lagrange equation}:
\begin{equation}\tag{EL${}^\alpha_h$}\label{elfh}
\forall k=1,\ldots,N-1, \; \dfrac{\partial L}{\partial x} \big( Q_k,(\DDM \Q)_k,t_k \big) + \DDP \left( \dfrac{\partial L}{\partial v} ( \Q,\DDM \Q, \T ) \right)_k = 0.
\end{equation}
\end{theorem}

\begin{proof}
Let $\Q \in (\R^d)^{N+1}$ and $\W \in (\R^d)^{N+1}_0$. Let us define the following function:
\begin{equation}
\fonction{\varphi}{\R}{\R}{\eps}{\LL^\alpha_h (\Q + \eps \W) = h \di \sum_{k=1}^N L \big( Q_k + \eps W_k, (\DDM \Q)_k + \eps (\DDM \W)_k,t_k \big).}
\end{equation}
Then, since $D\LL^\alpha_h (\Q)(\W) = \dot{\varphi} (0)$, we have:
\begin{equation}
D\LL^\alpha_h (\Q)(\W) = h \di \sum_{k=1}^N \left[ \dfrac{\partial L}{\partial x} \big( Q_k,(\DDM \Q)_k,t_k \big) \cdot W_k + \dfrac{\partial L}{\partial v} \big( Q_k,(\DDM \Q)_k,t_k \big) \cdot (\DDM \W)_k \right].
\end{equation}
Since $\W \in (\R^d)^{N+1}_0$ and using the discrete fractional integration by parts given in Lemma~\ref{lem2}, we obtain:
\begin{equation}
D\LL^\alpha_h (\Q)(\W) = h \di \sum_{k=1}^{N-1} \left[ \dfrac{\partial L}{\partial x} \big( Q_k,(\DDM \Q)_k,t_k \big) + \DDP \left( \dfrac{\partial L}{\partial v} ( \Q,\DDM \Q, \T ) \right)_k \right] \cdot W_k,
\end{equation}
which completes the proof.
\end{proof}

Theorem~\ref{thm1} completes the second step. Precisely, the discrete fractional Euler-Lagrange equation \eqref{elfh} is the variational integrator constructed for the fractional Euler-Lagrange equation \eqref{elf}. It is a numerical scheme preserving the variational structure of \eqref{elf} in the sense that the discrete solutions of \eqref{elfh} coincide with the discrete critical points of the discrete version $\LL^\alpha_h $ of $\LL^\alpha$.

\begin{remark}
We note that the discrete Euler-Lagrange equation \eqref{elfh} coincides with the numerical scheme \eqref{elfhpasbon} obtained with a direct discretization. It is important to note that such a phenomena is not obvious. Indeed, in the classical case $\alpha =1$, we have $D^\alpha_- = -D^\alpha_+ = d/dt$ and then, it is not obvious to know how to replace $d/dt$ at the discrete level. Indeed, one can choose $\Delta^1_-$ or $-\Delta^1_+$ or a mixing of the two of them. In the fractional case, the non locality of the fractional operators does not permit such a choice.
\end{remark}

In the classical case $\alpha = 1$, we note that \eqref{elfh} coincides with the classical discrete Euler-Lagrange equation provided in literature about variational integrators, see \cite{lubi,mars}.

\subsection{The fractional Dirichlet example}\label{section34}
In general, the fractional Euler-Lagrange equations are very difficult to solve explicitly. There are only few examples where exact solutions are provided and furthermore not in an explicit way. Due to the additional complexity of fractional operators, the purpose of this section is only to obtain experimental results and not to provide numerical analyses. \\

In this section, we consider the fractional Dirichlet example studied in \cite{agra} and where a \textit{quasi-explicit} solution is known. Precisely, we consider the interval $[a,b] = [0,1]$, $d=1$ and the following Lagrangian:
\begin{equation}
\fonction{L}{\R \times \R \times [0,1]}{\R}{(x,v,t)}{\dfrac{1}{2}v^2.}
\end{equation}
The associated fractional Euler-Lagrange equation \eqref{elf} is given by $D^\alpha_+ \circ D^\alpha_- q = 0$. For any $1/2 < \alpha \leq 1$, under boundary assumptions $q(0)=0$ and $q(1)=1$, this equation admits an unique solution given by:
\begin{equation}
q(t) = (2\alpha - 1) \di \int_0^t \dfrac{1}{[(1-x)(t-x)]^{1-\alpha}} \; dx.
\end{equation}
We note that for $\alpha =1$, the unique solution is $q(t) = t$. Nevertheless, in the strict fractional case $ 1/2 < \alpha < 1$, $q$ is not explicitly given since it is written with an improper integral. As a consequence, the following numerical results consider an approximation of the exact solution $q$ using high-order global adaptive quadrature. \\

Let us solve the discrete fractional Euler-Lagrange equation \eqref{elfh} given by:
\begin{equation}
\forall k=1,\ldots,N-1, \; (\Delta^\alpha_+ \circ \Delta^\alpha_- \Q)_k = 0,
\end{equation}
with boundary values $Q_0 = 0$ and $Q_N = 1$. For any $1/2 < \alpha \leq 1$ and any $ N \geq 2$, the $\ell^\infty$-error and the $\ell^2$-error between the exact solution $q$ and the discrete one $\Q$ are defined by:
\begin{equation}
E^\infty_{\alpha,N} = \di \max_{k=0,\ldots,N} \vert Q_k - q(t_k) \vert \; \text{and} \; E^2_{\alpha,N} = \sqrt{h \di \sum_{k=0}^N \vert Q_k - q(t_k) \vert^2 }, 
\end{equation}
where $h=1/N$. We compare $q$ and $\Q$ for $\alpha=3/4$ and $N=100$. The solutions are displayed on the following picture:
\begin{equation*}
\begin{array}{c}
\includegraphics[width=0.6\textwidth]{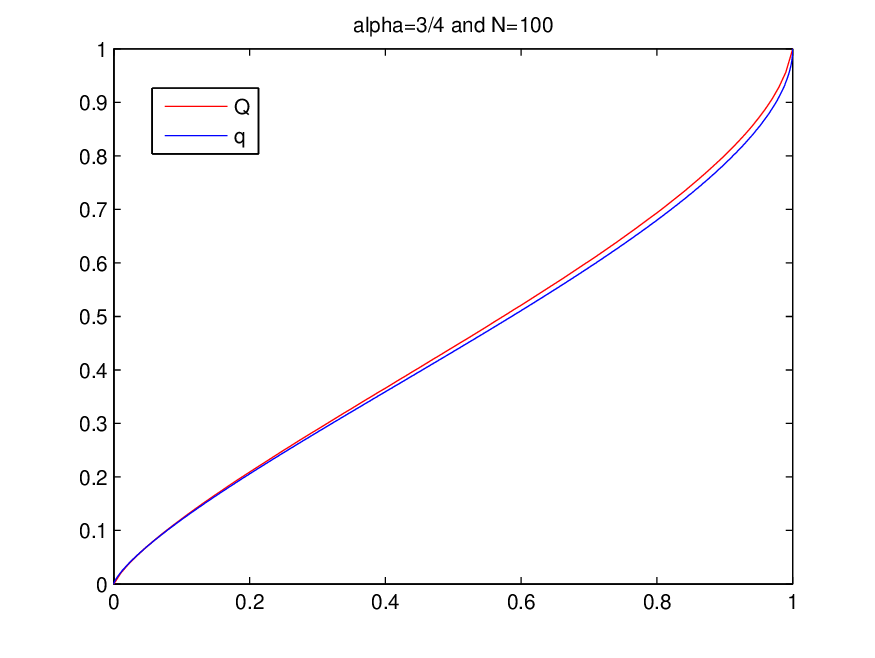}
\end{array}
\end{equation*}
The errors $E^\infty_{\alpha,N}$ and $E^2_{\alpha,N}$ for $\alpha =3/4$ and for varying $N$ are given in the following table:
\begin{center}
\begin{tabular}{|c|c|c|c|c|c|c|c|c|}
\hline $N$ & 50 & 100 & 200 & 250 & 500 & 1000 & 2000 & 4000  \\ \hline
$E^\infty_{3/4,N}$ & 0.0255 & 0.0185 & 0.0134 & 0.0120 & 0.0086 & 0.0062 & 0.0044 & 0.0031 \\ \hline
$E^2_{3/4,N}$ & 0.0140 & 0.0100 & 0.0072 & 0.0064 & 0.0046 & 0.0032 & 0.0023 & 0.0016 \\ \hline
\end{tabular}
\end{center}
Finally, the graphic representations of $\log(E^\infty_{3/4,N})$ and $\log(E^2_{3/4,N})$ with respect to $\log (h)$ are respectively given on the following pictures:
\begin{equation*}
\begin{array}{cc}
\includegraphics[width=0.4\textwidth]{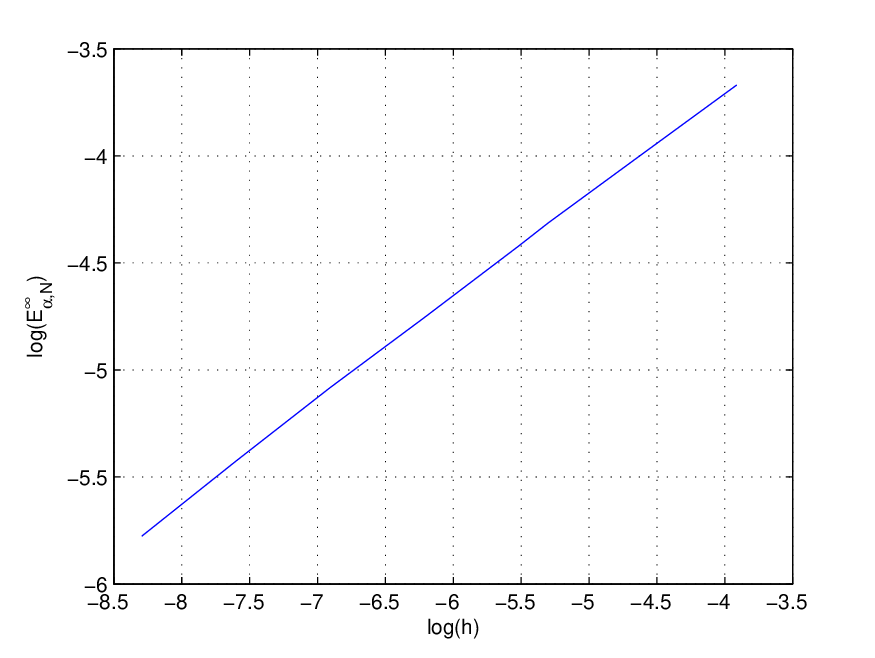} & \includegraphics[width=0.4\textwidth]{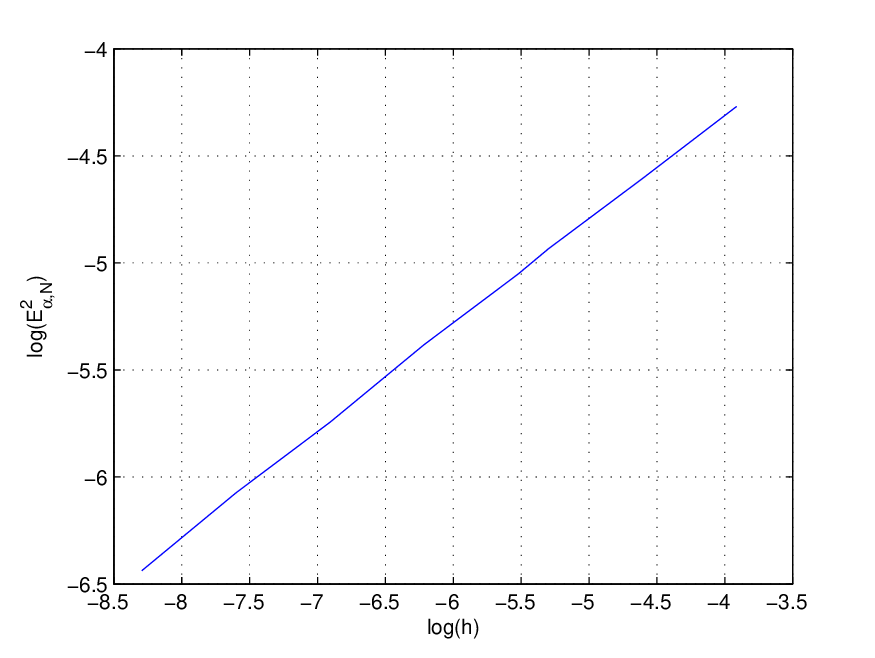}
\end{array}
\end{equation*}
We clearly obtain two linear functions with the common slope $ \lambda_{3/4} \simeq 0.51$. In this case, we conclude that we obtain an experimental convergence of order $ \lambda_{3/4} \simeq 0.51$. \\

Now, let us study the evolution of the errors $E^\infty_{\alpha,N}$ and $E^2_{\alpha,N}$ and the slope $\lambda_{\alpha}$ with respect to $\alpha$. Let us precise that the numerical study of the case $\alpha =1$ is not interesting because, in this case, the discrete Euler-Lagrange equation {\rm (EL${}^1_h$)} gives directly the exact solution $q(t) =t$. In the strict fractional case $1/2 < \alpha < 1$, the errors $E^\infty_{\alpha,N}$ and $E^2_{\alpha,N}$ are given in the following table:
\begin{center}
\begin{tabular}{|c|c|c|c|c|c|c|c|c|}
\hline  $\alpha \backslash N$ & 50 & 100 & 200 & 250 & 500 & 1000 & 2000 & 4000 \\ \hline
0.55 & 0.3811 & 0.3679 & 0.3527 & 0.3475 & 0.3310 & 0.3141 & 0.2974 & 0.2810 \\
 & 0.1746 & 0.1560 & 0.1393 & 0.1344 & 0.1203 & 0.1079 & 0.0972 & 0.0878 \\ \hline
0.6 & 0.1988 & 0.1788 & 0.1596 & 0.1537 & 0.1363 & 0.1204 & 0.1061 & 0.0933 \\
 & 0.0969 & 0.0825 & 0.0701 & 0.0665 & 0.0566 & 0.0483 & 0.0413 & 0.0354 \\ \hline
0.65 & 0.1023 & 0.0856 & 0.0711 & 0.0669 & 0.0552 & 0.0453 & 0.0372 & 0.0304 \\
 & 0.0525 & 0.0423 & 0.0340 & 0.0317 & 0.0255 & 0.0205 & 0.0165 & 0.0133 \\ \hline
0.7 & 0.0517 & 0.0403 & 0.0311 & 0.0286 & 0.0219 & 0.0168 & 0.0129 & 0.0098 \\
 & 0.0276 & 0.0210 & 0.0159 & 0.0146 & 0.0110 & 0.0083 & 0.0063 & 0.0048 \\ \hline
0.75 & 0.0255 & 0.0185 & 0.0134 & 0.0120 & 0.0086 & 0.0062 & 0.0044 & 0.0031 \\
 & 0.0140 & 0.0100 & 0.0072 & 0.0064 & 0.0046 & 0.0032 & 0.0023 & 0.0016 \\ \hline
0.8 & 0.0122 & 0.0083 & 0.0056 & 0.0049 & 0.0033 & 0.0022 & 0.0015 & 0.00097 \\
 & 0.0067 & 0.0045 & 0.0030 & 0.0027 & 0.0018 & 0.0012 & 0.0008 & 0.00052 \\ \hline
0.85 & 0.0055 & 0.0035 & 0.0022 & 0.0019 & 0.0012 & 0.00074 & 0.00046 & 0.00028 \\
 & 0.0030 & 0.0019 & 0.0012 & 0.0010 & 0.0006 & 0.0004 & 0.00025 & 0.00015 \\ \hline
0.9 & 0.0022 & 0.0013 & 0.00079 & 0.00066 & 0.00039 & 0.00022 & 0.00013 & 0.000075 \\
 & 0.0012 & 0.0007 & 0.00041 & 0.00035 & 0.0002 & 0.00012 & 0.00007 & 0.00004 \\ \hline
0.95 & 0.00068 & 0.00038 & 0.00021 & 0.00017 & 0.000095 & 0.000052 & 0.000028 & 0.000015 \\
 & 0.00037 & 0.0002 & 0.00011 & 0.00009 & 0.00005 & 0.000026 & 0.000014 & 0.000007 \\ \hline
\end{tabular}
\end{center}
We clearly see that the more $\alpha$ is close to $1$, the more the approximation is sensibly better. Now, we study the evolution of slope $\lambda_\alpha$ with respect to $\alpha$:
\begin{center}
\begin{tabular}{|c|c|c|c|c|c|c|c|c|c|}
\hline  $\alpha$ & 0.55 & 0.6 & 0.65 & 0.7 & 0.75 & 0.8 & 0.85 & 0.9 & 0.95 \\ \hline
$\lambda_{\alpha}$ & 0.11 & 0.2 & 0.3 & 0.39 & 0.51 & 0.62 & 0.72 & 0.8 & 0.91 \\ \hline
\end{tabular}
\end{center}
For this example, we conclude that we obtain an experimental convergence of order $ \lambda_\alpha = 2 \alpha - 1 $.

\section{Conservation of a fractional Noether-type result}\label{section4}
It is well-known that the conservation at the discrete level of the variational structure allows to preserve some properties relative to this structure. For example, in the classical case $\alpha =1$, we know that the classical Noether's theorem is preserved, see \cite{lubi,mars}. In this section, in the fractional case $0 < \alpha \leq 1$, we prove the conservation at the discrete level of the fractional Noether-type result proved simultaneously by Cresson \cite{cres6} and Torres \textit{et al.} \cite{torr3}.

\subsection{Reminder about a fractional Noether-type theorem}\label{section41}
We first review the definition of a one parameter group of diffeomorphisms: 
\begin{definition}
For any real $s$, let $\fonctionsansdef{\phi (s,\cdot)}{\R ^d}{\R ^d}$ be a diffeomorphism. Then, $\Phi = \{ \phi (s,\cdot) \}_{s \in \R}$ is a one parameter group of diffeomorphisms if it satisfies 
\begin{enumerate}
\item $\phi (0,\cdot) = Id_{\R ^d}$, 
\item $\forall s,s' \in \R, \; \phi (s,\cdot) \circ \phi (s',\cdot) = \phi (s+s',\cdot) $ ,
\item $\phi$ is of class $\mathcal{C}^2$.
\end{enumerate}
\end{definition}
Classical examples of one parameter groups of diffeomorphisms are given by translations and rotations.\\

The action of a one parameter group of diffeomorphisms on a Lagrangian allows to define the notion of a symmetry for a fractional Euler-Lagrange equation:

\begin{definition}
\label{definvariancefrac}
Let $\Phi = \{ \phi (s,\cdot) \}_{s \in \R}$ be a one parameter group of diffeomorphisms and let $L$ be a Lagrangian. $L$ is said to be $D^\alpha_-$-invariant under the action of $\Phi$ if it satisfies:
\begin{equation}
\forall q \text{ solution of \eqref{elf}},  \; \forall s \in \R, \; L \Big( \phi(s,q), D^\alpha_- \big( \phi (s,q) \big),t \Big) = L \big( q,D^\alpha_- q,t \big) .
\end{equation}
\end{definition}

For example, any quadratic Lagrangian is $D^\alpha_-$-invariant under the action of rotations. We refer to Section~\ref{section43} for the fractional harmonic oscillator example. We also note that the fractional Dirichlet example studied in Section~\ref{section34} is associated to a quadratic Lagrangian and then, it admits a symmetry. \\

With this definition of symmetry, Cresson \cite{cres6} and Torres \textit{et al.} \cite{torr3} proved the following result:

\begin{theorem}[Fractional Noether-type theorem]\label{lemcresson}
Let $L$ be a Lagrangian $D^\alpha_-$-invariant under the action of a one parameter group of diffeomorphisms $\Phi = \{ \phi (s,\cdot) \}_{s \in \R}$. Then, the following equality holds for any solution $q$ of \eqref{elf}:
\begin{equation}\label{eqlemcresson}
D^\alpha_- \left( \dfrac{\partial \phi}{\partial s} (0,q) \right) \cdot \dfrac{\partial L}{\partial v}(q,D^\alpha_- q,t) - \dfrac{\partial \phi}{\partial s} (0,q) \cdot D^\alpha_+ \left( \dfrac{\partial L}{\partial v}(q,D^\alpha_- q,t) \right)  = 0 .
\end{equation}
\end{theorem}

In the classical case $\alpha =1$, the classical Leibniz formula allows to rewrite \eqref{eqlemcresson} as the derivative of a product. Precisely, Theorem~\ref{lemcresson} leads to the classical Noether's theorem given by:

\begin{theorem}[Classical Noether's theorem]\label{thmnoetherclass}
Let $L$ be a Lagrangian $d/dt$-invariant under the action of a one parameter group of diffeomorphisms $\Phi = \{ \phi (s,\cdot) \}_{s \in \R}$. Then, the following equality holds for any solution $q$ of {\rm (EL${}^1$)}:
\begin{equation}\label{eqthmnoetherclass}
\dfrac{d}{dt} \left( \dfrac{\partial \phi}{\partial s} (0,q)  \cdot \dfrac{\partial L}{\partial v}(q,\dot{q},t) \right)  = 0, 
\end{equation}
where $\dot{q}$ is the classical derivative of $q$, \textit{i.e.} $dq/dt$.
\end{theorem}

Theorem~\ref{thmnoetherclass} provides an explicit constant of motion for any classical Euler-Lagrange equations {\rm (EL${}^1$)} admitting a symmetry. In the fractional case, such a simple Leibniz formula allowing to rewrite \eqref{eqlemcresson} as a total derivative with respect to $t$ is not known yet. Nevertheless, from Theorem~\ref{lemcresson}, we prove in \cite{bour2} the existence of an explicit constant of motion written with an infinite sum (this result is proved \textit{via} an iterative application of Leibniz formulas). For sake of brevity, we do not develop this result in this paper.

\subsection{A discrete fractional Noether-type theorem}
In this section, we prove that Theorem~\ref{lemcresson} is preserved at the discrete level with the discrete fractional Euler-Lagrange equations \eqref{elfh}. \\

As in the continuous case, a discrete symmetry of a discrete fractional Euler-Lagrange equation is based on the action of a one parameter group of transformations on the associated Lagrangian:

\begin{definition}
\label{definvariancedfrac}
Let $\Phi = \{ \phi (s,\cdot) \}_{s \in \R}$ be a one parameter group of diffeomorphisms and let $L$ be a Lagrangian. $L$ is said to be $\Delta^\alpha_-$-invariant under the action of $\Phi$ if it satisfies:
\begin{equation}\label{eqinvarfracd}
\forall \Q \text{ solution of \eqref{elfh}},  \; \forall s \in \R, \; L \Big( \phi(s,\Q), \Delta^\alpha_- \big( \phi (s,\Q) \big),\T \Big) = L \big( \Q,\Delta^\alpha_- \Q,\T \big) .
\end{equation}
\end{definition}

As in the continuous case, any quadratic Lagrangian is $\Delta^\alpha_-$-invariant under the action of rotations. We refer to Section~\ref{section43} for the fractional harmonic oscillator example. The fractional Dirichlet example studied in Section~\ref{section34} is also associated to a quadratic Lagrangian and then, it admits a discrete symmetry. \\

From this notion of discrete symmetry, we prove the discrete analogous version of Theorem~\ref{lemcresson}:

\begin{theorem}[Discrete fractional Noether-type theorem]\label{lemcressond}
Let $L$ be a Lagrangian $\Delta^\alpha_-$-invariant under the action of a one parameter group of diffeomorphisms $\Phi = \{ \phi (s,\cdot) \}_{s \in \R}$. Then, the following equality holds for any solution $\Q$ of \eqref{elfh}:
\begin{equation}\label{eqlemcressond}
\Delta^\alpha_- \left( \dfrac{\partial \phi}{\partial s} (0,\Q) \right) \cdot \dfrac{\partial L}{\partial v}(\Q,\Delta^\alpha_- \Q,\T) - \dfrac{\partial \phi}{\partial s} (0,\Q) \cdot \Delta^\alpha_+ \left( \dfrac{\partial L}{\partial v}(\Q,\Delta^\alpha_- \Q, \T) \right)  = 0 .
\end{equation}
\end{theorem}

\begin{proof} 
This proof is a direct adaptation to the discrete case of the proof of Theorem~\ref{lemcresson}. Let $\Q \in (\R^d)^{N+1}$ be a solution of \eqref{elfh}. Let us differentiate equation \eqref{eqinvarfracd} with respect to $s$ and invert the operators $\Delta^{\alpha}_-$ and $\partial / \partial s$. We finally obtain for any $s \in \R$ and any $k \in \{ 1,\ldots,N-1\}$:
\begin{multline}\label{astast}
\Delta^{\alpha}_- \left( \dfrac{\partial \phi}{\partial s} (s,\Q) \right)_k \cdot \dfrac{\partial L}{\partial v} \Big( \phi (s,Q_k), \Delta^\alpha_- \big( \phi (s,Q_k) \big),t_k \Big) \\ + \dfrac{\partial L}{\partial x} \Big( \phi (s,Q_k), \Delta^\alpha_- \big( \phi (s,Q_k) \big),t_k \Big) \cdot \dfrac{\partial \phi}{\partial s} (s,Q_k) = 0  .
\end{multline}
Since $\phi (0,\cdot) = Id_{\R^d}$ and $\Q$ is solution of \eqref{elfh}, taking $s=0$ in \eqref{astast} leads to \eqref{eqlemcressond}.
\end{proof}
From this last result, we conclude that the discrete fractional Euler-Lagrange equation \eqref{elfh} allows to preserve the fractional Noether's Theorem \ref{lemcresson} at the discrete level. \\

We recall the following discrete Leibniz formula:
\begin{equation}
\forall \boldsymbol{F}, \boldsymbol{G} \in \RN, \; \forall k=1,\ldots,N-1, \; \Delta^1_- \big( \boldsymbol{F}\cdot \sigma(\boldsymbol{G}) \big)_k = (\Delta^1_- \boldsymbol{F})_k \cdot G_k - F_k \cdot (\Delta^1_+ \boldsymbol{G})_k,
\end{equation}
where $\sigma (\boldsymbol{G})_k = G_{k+1}$ for any $k=0,\ldots,N-1$. From this discrete Leibniz formula, we note that Theorem~\ref{lemcressond} takes a particular simple expression in the classical case $\alpha =1$. It corresponds to the discrete version of Theorem~\ref{thmnoetherclass}:

\begin{theorem}[Discrete classical Noether's theorem]\label{thmnoetherdclass}
Let $L$ be a Lagrangian $\Delta^1_-$-invariant under the action of a one parameter group of diffeomorphisms $\Phi = \{ \phi (s,\cdot) \}_{s \in \R}$. Then, the following equality holds for any solution $\Q$ of {\rm (EL${}^1_h$)}:
\begin{equation}\label{dcnoether}
\Delta^1_- \left[ \dfrac{\partial \phi}{\partial s} (0,\Q) \cdot \sigma \left( \dfrac{\partial L}{\partial v} (\Q,\Delta^1_- \Q,\T ) \right) \right]= 0 .
\end{equation}
\end{theorem}
This result is a reformulation of the discrete Noether's theorem proved in \cite{lubi,mars}. We recall the following implication:
\begin{equation}\label{dconstant}
\forall \boldsymbol{F} \in \R^{N+1}, \; \Delta^1_- \boldsymbol{F} = 0 \Longrightarrow \exists c \in \R, \; \forall k=0,...,N, \; F_k = c.
\end{equation}
Consequently, as in the continuous case, Theorem~\ref{thmnoetherdclass} provides an explicit discrete constant of motion for any discrete classical Euler-Lagrange equations {\rm (EL${}^1_h$)} admitting a symmetry. In the fractional case, such a simple Leibniz formula allowing to rewrite \eqref{eqlemcressond} as a total discrete derivative is not known yet. Nevertheless, from Theorem~\ref{lemcressond}, we prove in \cite{bour2} the existence of a discrete constant of motion computable in a finite number of steps. For sake of brevity, we do not develop this result in this paper.

\subsection{The fractional harmonic oscillator example}\label{section43}
We conclude this section with the fractional harmonic oscillator example where Noether-type Theorems~\ref{lemcresson} and \ref{lemcressond} are available. Precisely, we consider $[a,b] = [0,1]$, $d=2$ and the following quadratic Lagrangian:
\begin{equation}
\fonction{L}{\R^2 \times \R^2 \times [0,1]}{\R}{(x,v,t)}{\dfrac{1}{2} (\Vert x \Vert^2 + \Vert v \Vert^2).}
\end{equation}
Since $L$ is quadratic, for any $0 < \alpha \leq 1$, $L$ is $D^\alpha_-$-invariant and $\Delta^\alpha_-$-invariant under the action of rotations given by:
\begin{equation}\label{rotation}
\fonction{\phi}{\R \times \R^2}{\R^2}{(s,x_1,x_2)}{\left( \begin{array}{cc} \cos (s) & - \sin (s ) \\ \sin (s) & \cos (s )  \end{array} \right) \left( \begin{array}{c} x_1 \\ x_2 \end{array} \right).}
\end{equation}

In the classical continuous case $\alpha =1$, the Euler-Lagrange equation {\rm (EL${}^1$)} is given by $\ddot{q}=q$ and consequently, the exact solutions are given by $q(t)=(q^1(t),q^2(t)) = (c_1 e^t+c_2e^{-t},c_3 e^t +c_4 e^{-t})$ for any $t \in [0,1]$ where $c_1$, $c_2$, $c_3$, $c_4 \in \R$. Applying Theorem~\ref{thmnoetherclass}, we obtain an explicit constant of motion. Indeed, for any $q = (q^1,q^2)$ solution of {\rm (EL${}^1$)}, there exists $c \in \R$ such that:
\begin{equation}
\dfrac{\partial \phi}{\partial s} (0,q)  \cdot \dfrac{\partial L}{\partial v}(q,\dot{q},t) = (-q^2,q^1) \cdot (\dot{q}^1,\dot{q}^2) =c.
\end{equation}
Similarly, the discrete Euler-Lagrange equation {\rm (EL${}^1_h$)} is given by $-\Delta^1_+ \circ \Delta^1_- \Q = \Q$ and Theorem~\ref{thmnoetherdclass} leads to an explicit discrete constant of motion. Indeed, for any $\Q = (\Q^1,\Q^2)$ discrete solution of {\rm (EL${}^1_h$)}, there exists $c \in \R$ such that:
\begin{equation}
\dfrac{\partial \phi}{\partial s} (0,\Q)  \cdot \sigma \left( \dfrac{\partial L}{\partial v}(\Q,\Delta^1_- \Q,\T) \right) = (-\Q^2,\Q^1) \cdot \sigma (\Delta^1_- \Q^1,\Delta^1_- \Q^2) =c.
\end{equation}

In the strict fractional continuous case $0 < \alpha < 1$, the fractional Euler-Lagrange equation \eqref{elf} is given by $q + \DP \circ \DM q =0$. Even if explicit solutions are not known yet, the fractional Noether-type Theorem~\ref{lemcresson} is available and leads to an explicit constant of motion in \cite{bour2}. In the discrete case, the discrete fractional Euler-Lagrange equation \eqref{elfh} is given by $\Q + \DDP \circ \DDM \Q =0$. The discrete solutions are computable and the discrete fractional Noether-type Theorem~\ref{lemcressond} (leading to an explicit discrete constant of motion in \cite{bour2}) is available.

\bibliographystyle{plain}


\begin{thebibliography}{10}

\bibitem{agra}
O.P. Agrawal.
\newblock Formulation of {E}uler-{L}agrange equations for fractional
  variational problems.
\newblock {\em J. Math. Anal. Appl.}, 272(1):368--379, 2002.

\bibitem{agra2}
O.P. Agrawal.
\newblock A general formulation and solution scheme for fractional optimal
  control problems.
\newblock {\em Nonlinear Dynam.}, 38(1-4):323--337, 2004.

\bibitem{agra3}
O.P. Agrawal.
\newblock A formulation and numerical scheme for fractional optimal control
  problems.
\newblock {\em J. Vib. Control}, 14(9-10):1291--1299, 2008.

\bibitem{agra4}
O.P. Agrawal.
\newblock A general finite element formulation for fractional variational
  problems.
\newblock {\em J. Math. Anal. Appl.}, 337(1):1--12, 2008.

\bibitem{alme}
R.~Almeida, A.B. Malinowska, and D.F.M. Torres.
\newblock A fractional calculus of variations for multiple integrals with
  application to vibrating string.
\newblock {\em J. Math. Phys.}, 51(3):033503, 12, 2010.

\bibitem{bagl}
R.~L. Bagley and R.~A. Calico.
\newblock Fractional order state equations for the control of viscoelastically
  damped structures.
\newblock {\em Journal of Guidance, Control, and Dynamics}, 14:304--311, 1991.

\bibitem{bale2}
D.~Baleanu and S.I. Muslih.
\newblock Lagrangian formulation of classical fields within
  {R}iemann-{L}iouville fractional derivatives.
\newblock {\em Phys. Scripta}, 72(2-3):119--121, 2005.

\bibitem{bour2}
L.~Bourdin, J.~Cresson, and I.~Greff.
\newblock A continuous/discrete fractional {N}oether's theorem.
\newblock {\em Communications in Nonlinear Science and Numerical Simulation},
  18(4):878 -- 887, 2013.

\bibitem{comt}
F.~Comte.
\newblock Op\'erateurs fractionnaires en \'econom\'etrie et en finance.
\newblock {\em Pr\'epublication MAP5}, 2001.

\bibitem{cres6}
J.~Cresson.
\newblock Fractional embedding of differential operators and {L}agrangian
  systems.
\newblock {\em J. Math. Phys.}, 48(3):033504, 34, 2007.

\bibitem{deft}
O.~Defterli.
\newblock A numerical scheme for two-dimensional optimal control problems with
  memory effect.
\newblock {\em Comput. Math. Appl.}, 59(5):1630--1636, 2010.

\bibitem{diet}
K.~Diethelm.
\newblock {\em The analysis of fractional differential equations}, volume 2004
  of {\em Lecture Notes in Mathematics}.
\newblock Springer-Verlag, Berlin, 2010.
\newblock An application-oriented exposition using differential operators of
  Caputo type.

\bibitem{dubo}
F.~Dubois, A-C. Galucio, and N.~Point.
\newblock Introduction \`a la d\'erivation fractionnaire. {T}h\'eorie et
  applications.
\newblock {\em S{\'e}rie des Techniques de l'ing{\'e}nieur}, 2009.

\bibitem{torr3}
G.S.F. Frederico and D.F.M. Torres.
\newblock A formulation of {N}oether's theorem for fractional problems of the
  calculus of variations.
\newblock {\em J. Math. Anal. Appl.}, 334(2):834--846, 2007.

\bibitem{lubi}
E.~Hairer, C.~Lubich, and G.~Wanner.
\newblock {\em Geometric numerical integration}, volume~31 of {\em Springer
  Series in Computational Mathematics}.
\newblock Springer-Verlag, Berlin, second edition, 2006.
\newblock Structure-preserving algorithms for ordinary differential equations.

\bibitem{hilf3}
R.~Hilfer.
\newblock Applications of fractional calculus in physics.
\newblock {\em World Scientific, River Edge, New Jersey}, 2000.

\bibitem{kilb}
A.A. Kilbas, H.M. Srivastava, and J.J. Trujillo.
\newblock {\em Theory and applications of fractional differential equations},
  volume 204 of {\em North-Holland Mathematics Studies}.
\newblock Elsevier Science B.V., Amsterdam, 2006.

\bibitem{liu}
Yanqin Liu and Baogui Xin.
\newblock Numerical solutions of a fractional predator-prey system.
\newblock {\em Adv. Difference Equ.}, pages Art. ID 190475, 11, 2011.

\bibitem{mach}
J.~Tenreiro Machado, Virginia Kiryakova, and Francesco Mainardi.
\newblock Recent history of fractional calculus.
\newblock {\em Commun. Nonlinear Sci. Numer. Simul.}, 16(3):1140--1153, 2011.

\bibitem{magi}
R.L. Magin.
\newblock Fractional calculus models of complex dynamics in biological tissues.
\newblock {\em Comput. Math. Appl.}, 59(5):1586--1593, 2010.

\bibitem{mars}
J.E. Marsden and M.~West.
\newblock Discrete mechanics and variational integrators.
\newblock {\em Acta Numer.}, 10:357--514, 2001.

\bibitem{meer}
M.M. Meerschaert and C.~Tadjeran.
\newblock Finite difference approximations for fractional advection-dispersion
  flow equations.
\newblock {\em J. Comput. Appl. Math.}, 172(1):65--77, 2004.

\bibitem{podl}
I.~Podlubny.
\newblock {\em Fractional differential equations}, volume 198 of {\em
  Mathematics in Science and Engineering}.
\newblock Academic Press Inc., San Diego, CA, 1999.
\newblock An introduction to fractional derivatives, fractional differential
  equations, to methods of their solution and some of their applications.

\bibitem{riew}
F.~Riewe.
\newblock Nonconservative {L}agrangian and {H}amiltonian mechanics.
\newblock {\em Phys. Rev. E (3)}, 53(2):1890--1899, 1996.

\bibitem{samk}
S.G. Samko, A.A. Kilbas, and O.I. Marichev.
\newblock {\em Fractional integrals and derivatives}.
\newblock Gordon and Breach Science Publishers, Yverdon, 1993.
\newblock Theory and applications, Translated from the 1987 Russian original.

\bibitem{neel}
A.~Zoia, M.-C. N\'eel, and A.~Cortis.
\newblock Continuous-time random-walk model of transport in variably saturated
  heterogeneous porous media.
\newblock {\em Phys. Rev. E}, 81(3):031104, Mar 2010.

\end{thebibliography}

\end{document}